\documentclass[12pt]{amsart}

\newtheorem{thm}{Theorem}[section]   
\newtheorem{cor}[thm]{Corollary}     
\newtheorem{lem}[thm]{Lemma}         
\newtheorem{prop}[thm]{Proposition}  

\theoremstyle{definition}
\newtheorem{defn}[thm]{Definition}   

\theoremstyle{remark}
\newtheorem{rem}[thm]{Remark}        
\newtheorem{ex}[thm]{Example}        

\numberwithin{equation}{section}     

\newcommand{\secref}[1]{Section~\ref{#1}}
\newcommand{\thmref}[1]{Theorem~\ref{#1}}

\newcommand{\lemref}[1]{Lemma~\ref{#1}}
\newcommand{\propref}[1]{Proposition~\ref{#1}}

\newcommand{\exref}[1]{Example~\ref{#1}}

\renewcommand{\H}{\mathcal H}
\newcommand{\B}{\mathcal B}
\newcommand{\K}{\mathcal K}
\renewcommand{\L}{\mathcal L}
\newcommand{\Z}{\mathbb Z}
\newcommand{\F}{\mathbb F}
\newcommand{\TO}{\mathcal T\!\mathcal O}
\renewcommand{\O}{\mathcal O}
\newcommand{\W}{\mathcal W}
\newcommand{\D}{\mathcal D}
\newcommand{\T}{\mathbb T}
\newcommand{\cstar}{$C^*$-}
\renewcommand{\star}{\text{${}^*$}-}
\newcommand{\spacetext}[1]{\quad\text{#1}\quad}
\newcommand{\righttext}[1]{\qquad\text{#1 }}
\newcommand{\noqed}{\renewcommand{\qed}{}}
\newcommand{\ad}{\operatorname{Ad}}
\newcommand{\clsp}{\overline{\text{sp}}}
\newcommand{\inv}{^{-1}}

\begin{document}

\title{Characterisations of crossed products by partial actions}
\author{John Quigg}
\address{Department of Mathematics\\Arizona State University\\
 Tempe, Arizona 85287}
\email{quigg@math.la.asu.edu}
\author{Iain Raeburn}
\address{Department of Mathematics\\University of Newcastle\\
 Newcastle, New South Wales 2308\\Australia}
\email{iain@frey.newcastle.edu.au}
\thanks{This research was partially supported by the Australian Research 
Council and by the National Science 
Foundation (under Grant No. DMS9401253).} 
\date{February 13, 1996}

\maketitle

\begin{abstract}
Partial actions of discrete groups on $C^*$-algebras and the associated
crossed products have been studied by Exel and McClanahan. We
characterise these crossed products in terms of the spectral subspaces
of the dual coaction, generalising and simplifying a theorem of Exel
for single partial automorphisms. We then use this characterisation to
identify the Cuntz algebras and the Toeplitz algebras of Nica as crossed
products by partial actions.
\end{abstract}

\section*{Introduction}

Exel has recently introduced and studied partial automorphisms of a
$C^*$-algebra $A$: isomorphisms of one ideal in $A$ onto another
\cite{exe:partial}. He has shown that many interesting $C^*$-algebras
can be viewed as crossed products by partial automorphisms, and that
these crossed products have much in common with ordinary crossed
products by actions of $\Z$. McClanahan subsequently extended Exel's
ideas to cover \emph{partial actions} of more general groups by partial
automorphisms, and showed that, rather surprisingly, many important
results on crossed products by free groups carry over to crossed
products by partial actions \cite{mcc}. 

Here we give a characterisation of (reduced) crossed products by partial
actions of discrete groups, which is similar in spirit to that of Landstad
for ordinary crossed products (see \cite{lan:dual} or \cite[7.8.8]{ped}),
and which both generalises and simplifies Exel's characterisation of
crossed products by single partial automorphisms \cite[Theorem
4.21]{exe:partial}. Our main result says that a $C^*$-algebra $B$ is a
crossed product by a partial action of $G$ if and only if it carries a
coaction $\delta$ of $G$ and there is a \emph{partial representation} of
$G$ by partial isometries in the double dual $B^{**}$ which induces
suitable isomorphisms among the spectral subspaces of $\delta$; this
result takes a particularly elegant form when $G$ is the free group
$\F_n$. We then use our classification to identify the Cuntz
algebras $\O_n$, the Cuntz-Krieger algebras $\O_A$, and the Toeplitz or
Wiener-Hopf algebras of Nica
\cite{nic} as crossed products by partial actions of $\F_n$. We have not
previously seen the canonical coaction of $\F_n$ on $\O_n$ used in a
serious way, so these ideas may have interesting implications for the
study of coactions of discrete groups. 

\medskip

We begin with a discussion of partial actions and covariant
representations. A partial action $\alpha$ of $G$ on $A$ is a collection
$\{D_s:s\in G\}$ of ideals in $A$, and isomorphisms $\alpha_s$ of
$D_{s^{-1}}$ onto $D_s$, such that $\alpha_{st}$ extends
$\alpha_s\circ\alpha_t$ from its natural domain
$D_{t^{-1}}\cap\alpha_t^{-1}(D_{s^{-1}})$. Calculations involving the
domains can be tricky, so we have taken care to make the various
relationships explicit. A partial representation of $G$ is a map $u$ of
$G$ into the set of partial isometries on a Hilbert space (or in a
$C^*$-algebra) such that $u_{st}$ extends $u_su_t$, and a covariant
representation
$(\pi,u)$ of $(A,G,\alpha)$ consists of a representation $\pi$ of $A$ and
a partial representation $u$ of $G$ such that
$\pi(\alpha_s(a))=u_s\pi(a)u_s^*$ for
$a\in D_{s^{-1}}$. McClanahan did not discuss partial representations in
their own right, so we have included a detailed discussion of them and
their relationship to covariant representations.

A key technical innovation in our treatment is the implementation of
Hilbert-module isomorphisms of spectral subspaces by multipliers of
imprimitivity bimodules, as introduced in \cite{ech-rae}; the particular
multipliers involved here will form the partial representation of $G$ in
the double dual of the crossed product. We therefore
recall in \S2 some facts about multipliers of bimodules, relate them to
Hilbert-module isomorphisms, and discuss how in certain situations the
whole structure can be embedded in the double dual of a $C^*$-algebra.

In \S3, we construct the crossed product $A\times_\alpha G$ of a
partial action $\alpha$, as the  $C^*$-algebra generated by a universal
covariant representation of $(A,G,\alpha)$ in $(A\times_\alpha
G)^{**}$. Associated to any faithful representation $\pi$ of $A$ is a
regular representation of $A\times_\alpha G$; up to isomorphism, the
image is independent of the choice of $\pi$, and is called the reduced
crossed product $A\times_{\alpha,r}G$. Of course, both crossed products
turn out to be the ones studied in \cite{mcc}, but our emphasis on
universal properties allows us to see quickly that they carry a dual
coaction of $G$. Our characterisation of the reduced crossed product in
terms of this dual coaction is Theorem \ref{duality}.

An ordinary action $\alpha$ of the free group $\F_n$ is determined
completely by the $n$ automorphisms $\alpha_{g_i}$ corresponding to
generators  $\{g_i\}$ of $\F_n$. It is quite easy to construct
partial actions of $\F_n$ from $n$ partial automorphisms \cite[Example
2.3]{mcc}, but in general even partial actions of $\F_1=\Z$ need not arise
this way. So we concentrate in \S5 on a family of partial actions $\alpha$
of $\F_n$ which are determined by $\{\alpha_{g_i}\}$; crossed products
by such \emph{multiplicative} partial actions can be characterised in
terms of the spectral subspaces of the dual coaction corresponding to the
generators of $\F_n$. The main result here is Theorem
\ref{elementaryduality}, and its applications to Cuntz algebras,
Cuntz-Krieger algebras and Nica's Toeplitz algebras are the content of our
last section.

\medskip

\noindent\emph{Acknowledgements.} As this paper was being written up, the
authors received a copy of
\cite{exe:twist}, in which Exel proves a result related to our
\thmref{duality}. However, his result concerns \cstar algebraic
bundles, and uses techniques substantially different from ours.

The bulk of this research was carried out while the first 
author was visiting the University of Newcastle. He wishes to thank
Iain  Raeburn for his hospitality. The authors thank Marcelo Laca for
helpful conversations, and  for introducing them to the
work of Nica through a collaboration with the second author \cite{lr}.  

\section{Partial actions and covariant representations}
\label{setup}

\begin{defn}
 A \emph{partial 
action} of a discrete group $G$ on a $C^*$-algebra $A$ consists of a
collection
$\{D_s\}_{s\in G}$ of  closed ideals of $A$ and 
isomorphisms
$\alpha_s\colon  D_{s\inv}\to D_s$ such that 
\begin{enumerate}
\item $D_e=A$;
\item $\alpha_{st}$ extends $\alpha_s\alpha_t$ for all $s,t\in G$ (where 
the domain of $\alpha_s\alpha_t$ is $\alpha\inv_t(D_{s\inv})$, which
is by definition contained in $D_{t^{-1}}$). 
\end{enumerate}
\end{defn}

We suppose for the rest of this section that $\alpha$ is a partial
action of $G$ on $A$. We shall frequently need to intersect domains of the
partial isomorphims, and shall use without comment that the intersection
of two ideals $I,J$ in a
$C^*$-algebra is the ideal $IJ:=\overline{\operatorname{sp}}\{ij:i\in
I,j\in J\}$. We shall also need to recall that the double dual $I^{**}$
of an ideal $I$ naturally embeds as an ideal in $A^{**}$.

\begin{lem}\label{domains}
For $s,t\in G$ we have:
\begin{enumerate}
\item $\alpha_e=\iota$, and $\alpha_{s^{-1}}=\alpha_s^{-1}$;
\item $\alpha_s(D_{s^{-1}}D_t)=D_sD_{st}$;
\item $\alpha_s\circ\alpha_t$ is an isomorphism of
$D_{t^{-1}}D_{t^{-1}s^{-1}}$ onto $D_sD_{st}$.
\end{enumerate}
\end{lem}

\begin{rem}
Suppose instead of (ii) we had
\[ \alpha_s(D_{s\inv}D_t)\subset D_{st}\quad\text{and}\quad
\alpha_{st}(x)=\alpha_s\alpha_t(x)\text{ for }
x\in D_{t\inv}D_{t\inv s\inv}, \]
so that $\alpha$ is a partial action in the sense of McClanahan. Then
$\alpha$ would satisfy (ii), and hence be a partial action in our sense.
To see this, just note that 
\[
\operatorname{dom}(\alpha_s\circ\alpha_t)=\alpha_t^{-1}(D_{s^{-1}})=
\alpha_t^{-1}(D_tD_{s^{-1}})\subset D_{t^{-1}}D_{t^{-1}s^{-1}},
\]
so the equation $\alpha_{st}=\alpha_s\circ\alpha_t$ on
$D_{t^{-1}}D_{t^{-1}s^{-1}}$ says  what we need. Indeed, by part (ii) of
the lemma, this says {\em exactly} what we need, so our definition is
equivalent to McClanahan's. 
\end{rem}

\begin{proof}[Proof of Lemma \ref{domains}]
For part (i), observe that $\alpha_e\alpha_e$ is defined everywhere and
equal to $\alpha_e$, which is therefore the identity transformation on
all of $A$. Because $\alpha_e$ extends $\alpha_s\alpha_{s^{-1}}$ and
$\alpha_{s^{-1}}\alpha_s$, this forces
$\alpha_{s^{-1}}=\alpha_s^{-1}:D_s=D_{(s^{-1})^{-1}}\to D_{s^{-1}}$. For
(ii), we note that $\alpha_{t^{-1}}(D_{s^{-1}})$ is by definition
$\alpha_{t^{-1}}(D_tD_{s^{-1}})$, and because $\alpha_{st}$ extends
$\alpha_s\alpha_t$, we have 
\begin{equation}\label{incldomains}
\alpha_{t^{-1}}(D_tD_{s^{-1}})\subset D_{t^{-1}}D_{t^{-1}s^{-1}}\ 
\text{ for all }s,t.
\end{equation}
Applying this with $t^{-1}$ replaced by $t$ shows
\[
\alpha_t\big(D_{t^{-1}}D_{t^{-1}s^{-1}}\big)\subset
D_tD_{t(t^{-1}s^{-1})}=D_tD_{s^{-1}},
\]
and because $\alpha_{t^{-1}}=\alpha_t^{-1}$, this implies that we must
have equality in (\ref{incldomains}). Since the left-hand side of
(\ref{incldomains}) is the natural domain of $\alpha_s\alpha_t$,
and the range of $\alpha_s\alpha_t$ is the natural domain of
$\alpha_t^{-1}\alpha_s^{-1}=\alpha_{t^{-1}}\alpha_{s^{-1}}$, part (iii)
follows.
\end{proof}

\begin{defn}
For $s\in G$, we let $p_s$ denote the projection
in $A^{**}$ which is the identity of $D^{**}_s$.
\end{defn}

 The projections $p_s$ belong to the center
of $A^{**}$, and $p_s$ is the weak* limit of any bounded approximate
identity for $D_s$. We always have $p_s\in M(D_s)$, but $p_s$ may not be
in
$M(A)$, as shown by the following example.

\begin{ex}
\label{shift}
Let $A=C_0(0,\infty)$ and $G=\Z_2$. Define $D_1 = \{ f \in A \mid f(x) = 
0 \text{ for } x\le 1 \}$, and let $\alpha_1$ be the identity map of $D_1$.
Then $p_1$ is the characteristic function of $(1,\infty)$, which is not in 
$M(A)$ because it is not continuous on $(0,\infty)$. 
\end{ex}

To define covariant representations of partial actions, we need
an appropriate notion of \emph{partial representations} of
groups by partial isometries.  The idea is  that $u_{st}$ should
extend $u_su_t$; for this to make sense, $u_su_t$ must be a partial
isometry, so we insist that the range projections commute. The following
Lemma describes what we mean by ``$v$ extends $u$'', and is presumably
standard.

\begin{lem} We define a relation $\preceq$ on the set of partial
isometries on a Hilbert space $\H$ by
\[
u\preceq v\Longleftrightarrow uu^*=uv^*.
\]
Then $u\preceq v$ precisely when the initial space $u^*u(\H)$ of $u$ is
contained in $v^*v(\H)$, and $v=u$ on $u^*u(\H)$; we have $u\preceq
v\Longleftrightarrow u^*u=v^*u$, and $\preceq$ is a partial order on the
set of partial isometries on $\H$.
\end{lem}

\begin{proof}
Since $uu^*$ is self-adjoint, $uu^*=uv^*$ implies $uu^*=vu^*$. But then
$vu^*u=(uu^*)u=u$ implies that $v=u$ on the range of $u^*u$. In
particular, this implies that $v^*$ maps the range of $u$ into the
initial space of $u$, and hence $uu^*=uv^*$ implies
\[
u^*u=u^*(uu^*)u=u^*(uv^*)u=(u^*u)v^*u=v^*u.
\]
Conversely, $u^*u=v^*u$ implies that $v^*=u^*$ on the range of $uu^*$,
and hence that $v$ maps the initial space of $u$ into the range of $u$;
thus
\[
uu^*=u(u^*u)u^*=u(u^*v)u^*=(uu^*)vu^*=vu^*.
\]
Finally, since $u\preceq v$ implies $uv^*v=u$, it is easy to check that
$\preceq$ is transitive, and that $u\preceq v$, $v\preceq u$ force $u=v$.
\end{proof}

\begin{defn}
\label{partrep}
A \emph{partial representation} of $G$ on a Hilbert space $\H$ is a
map $u\colon G\to\B(\H)$ such that the $u_s$ are partial isometries
with commuting range projections, and
\begin{alignat}{2}
\label{identity}
u_eu^*_e  & = 1;\\
\label{inverse}
u^*_su_s & = u_{s\inv}u^*_{s\inv} && \righttext{for all} s \in G; \\
\label{hom}
u_su_t & \preceq u_{st} && \righttext{for all} s,t \in G.
\end{alignat}
\end{defn}

We begin by listing some straightforward consequences of the definition. 

\begin{lem}
\label{altpartrep}
If $u$ is a partial representation, then
\begin{alignat}{2}
\label{altidentity}
u_e & = 1; \\
\label{altinverse}
u^*_s & = u_{s\inv} && \righttext{for all} s \in G; \\
\label{althom}
u_su_t & = u_su^*_su_{st} && \righttext{for all} s,t \in G.
\end{alignat}
\end{lem}

\begin{proof}For \eqref{altidentity}, note that \eqref{identity} 
and \eqref{hom} imply that $u_e$ is an idempotent coisometry. 
Since $u^*_s$ and $u_{s\inv}$ are partial isometries with the same range 
projection, and $u_su_{s\inv} \preceq 1$, we have \eqref{altinverse}. 
For the last part, we use the relation $u\preceq v$ in the forms
$u=uu^*v$, $uu^*=uv^*$, and then $v^*u=u^*u$, to deduce that
\begin{align*}
u_su_t
& = (u_su_t)(u_su_t)^*u_{st}= (u_su_tu^*_tu^*_s)(u_su^*_su_{st})\\
&=(u_su_tu_{st}^*)(u_su_s^*u_{st})=u_s(u_tu_{st}^*u_s)u_s^*u_{st}\\
&=u_s(u_s^*u_{st}u_{st}^*u_s)u_s^*u_{st},
\end{align*}
which equals $u_su_s^*u_{st}$ because $u_s^*u_{st}$ is a partial
isometry.
\end{proof}

\begin{rem}
Conditions \eqref{altidentity}--\eqref{althom} are stronger forms of
\eqref{identity}--\eqref{hom}; to see this for \eqref{hom}, note that
$u_su^*_s$ is a projection commuting with $u_{st}u^*_{st}$, and hence
\[
(u_su_t)u_{st}^*=u_su_s^*u_{st}u_{st}^*=
(u_su_s^*u_{st})(u_{st}^*u_su_s^*)=(u_su_t)(u_su_t)^*.
\]
\end{rem}

\begin{defn}
\label{covrep}
Let $\alpha$ be a partial action of $G$ on $A$. A \emph{covariant
representation} of $(A,G,\alpha)$ is a pair $(\pi,u)$ consisting of a
nondegenerate representation $\pi$ of $A$ and a partial representation
$u$ of $G$ on the same Hilbert space, satisfying
\begin{gather}
\label{rangeproj}
u_su^*_s = \pi(p_s); \\
\label{covariance}
\pi(\alpha_s(a))=u_s\pi(a)u_s^*
\ \text{ for } a \in D_{s\inv}.
\end{gather}
\end{defn}

\begin{lem}
\label{altcovrep}
Let $\alpha$ be a partial action of $G$ on $A$, let $\pi \colon A \to 
\B(\H)$ be a nondegenerate representation, and let $u \colon G \to \B(\H)$ 
be a map satisfying \eqref{rangeproj}--\eqref{covariance}. Then the 
following are equivalent\textup: 
\begin{enumerate}
\item $u$ is a partial
representation \textup(and $(\pi,u)$ is a covariant 
representation\textup{);}
\item $u_su_t \preceq u_{st}$ for all $s,t \in G$\textup;
\item $\pi(p_{st})u_su_t = \pi(p_s)u_{st}$ for all $s,t \in G$\textup;
\item $\pi(a)u_su_t = \pi(a)u_{st}$ for all $s,t \in G,a \in D_sD_{st}$.
\end{enumerate}
\end{lem}

\begin{proof}
Note first that by \eqref{rangeproj}--\eqref{covariance}, the 
$u_s$ are partial isometries with commuting range projections, and 
$u_eu^*_e = \pi(p_e) = 1$, so \eqref{identity} holds. Further, (iii) is 
equivalent to (iv) since $p_sp_{st}$ is the identity of $(D_sD_{st})^{**}$ 
in $A^{**}$, and \lemref{altpartrep} tells us that (i) implies (iii). 

Assume (ii). We will show \eqref{altinverse}, giving \eqref{inverse}, hence 
(i). First note that \eqref{altidentity} holds, since its proof only 
required \eqref{identity} and \eqref{hom}. We have 
\begin{equation*}
\begin{split}
u_su^*_s
& = \pi(p_s)
= \pi\circ\alpha_s(p_{s\inv})
= \ad u_s\circ\pi(p_{s\inv}) \\
& = u_su_{s\inv}u^*_{s\inv}u^*_s
= u_su_{s\inv}u^*_{ss\inv}
= u_su_{s\inv},
\end{split}
\end{equation*}
so that
\begin{equation}
\label{less}
u_s \preceq u^*_{s\inv}.
\end{equation}
Since the partial ordering $\preceq$ on partial isometries is
conjugation-invariant,  we get $u^*_s \preceq u_{s\inv}$. Applying 
\eqref{less} with $s$ replaced by $s\inv$, we arrive at
\[ u^*_s \preceq u_{s\inv} \preceq u^*_s, \]
so $u^*_s = u_{s\inv}$, which is \eqref{altinverse}.

Finally, assume (iii). To show (ii), we again need \eqref{altidentity} and 
\eqref{altinverse}. \eqref{altidentity} follows from \eqref{identity} and 
(iii). For 
\eqref{altinverse}, we have 
\[ u_su_{s\inv}
= \pi(p_{ss\inv})u_su_{s\inv}
= \pi(p_s)u_{ss\inv}
= u_su^*_s, \]
and the argument of the preceding paragraph shows $u_{s^{-1}}=u_s^*$. 
We now show (ii):
\begin{equation*}
\begin{split}
u_su_tu^*_tu^*_s
& = u_s\pi(p_{s\inv})\pi(p_t)u^*_s
= \pi \circ \alpha_s(p_{s\inv}p_t)
= \pi(p_sp_{st}) \\
& = \pi(p_sp_{st})u_{st}u^*_{st}
= \pi(p_sp_{st})u_su_tu^*_{st} \\
& = u_s\pi(p_{s\inv}p_t)u^*_su_su_tu^*_{st}
= u_su_tu^*_{st}. \qed
\end{split}
\end{equation*}
\noqed
\end{proof}

\begin{rem}
The previous Lemma shows that our definition of covariant
representation is equivalent to McClanahan's \cite{mcc}. So ours
is actually a slight improvement over McClanahan's, in that the
conditions 
$u^*_su_s = \pi(p_{s\inv})$ and $u^*_s = u_{s\inv}$ follow automatically.
\end{rem}

\section{Multipliers of imprimitivity bimodules}\label{mbimods}

Recall from \cite{ech-rae} that if $X$ is a
$C - D$ imprimitivity bimodule, a \emph{multiplier} of $X$ is a 
pair $m = (m_C,m_D)$, where $m_C \in \L_C(C, X)$ and
$m_D \in \L_D(D,X)$ satisfy 
\[ m_C(c)\cdot d = c \cdot  m_D(d) \spacetext{for} c \in C,d \in D. \]
(Actually, \cite[Lemma 1.4]{ech-rae} shows that adjointability of $m_C$
and 
$m_D$ is automatic.) The set $M(X)$ of multipliers of $X$ is 
called the \emph{multiplier bimodule}; with
\[ c \cdot m = m_C(c) \spacetext{and} m\cdot d = m_D(d), \]
$M(X)$ becomes a $C-D$ bimodule containing $X$. The module 
actions of $C$ and $D$ on $M(X)$ extend to $M(C)$ and
$M(D)$, and the $C$- and $D$-valued inner products on $X$ 
extend to $M(C)$- and $M(D)$-valued inner products on $M(X)$, which we
continue to denote by
${}_C\langle \cdot,\cdot \rangle$ and $\langle \cdot,\cdot \rangle_D$. 
For $x \in X$ and $m \in M(X)$ we have 
\begin{equation}\label{ipsonm(x)}
 {}_C\langle x, m \rangle = m^*_C(x)
\spacetext{and}
\langle m, x \rangle_D = m^*_D(x). 
\end{equation}

\begin{lem}\label{module iso}
Let $X$ be  a $C-D$ imprimitivity bimodule.
\begin{enumerate}
\item
There is a left Hilbert $C$-module isomorphism $\psi$ of $X$ onto $C$ if
and only if there is a multiplier $m=(m_C,m_D)$ of $X$ such that
\begin{equation}\label{mult}
{}_C\langle m,m\rangle=1_{M(C)},\ \ \langle m,m\rangle_D=1_{M(D)}
\end{equation}
and $\psi(x)={}_C\langle x,m\rangle=m_C^*(x)$ for all $x\in X$.
\item
Let $m$ be a multiplier of $X$ satisfying (\ref{mult}), and let
$\psi=m_C^*:X\to C$ and $\phi=m_D^*:X\to D$ be the corresponding
isomorphisms of Hilbert modules. Then there is a $C^*$-algebra isomorphism
$\alpha$ of
$D$ onto $C$ such that $\alpha(d)={}_C\langle m\cdot d,m\rangle$ and
$\psi=\alpha\circ \phi$. We write $\alpha=\ad m$.
\end{enumerate}
\end{lem}

\begin{proof}
Since $\psi$ preserves the $C$-valued inner products, the inverse
$\psi^{-1}:C\to X$ is an adjoint for $\psi$, and $\psi^{-1}$ itself is
in $\L_C(C,X)$. Thus by \cite[Proposition 1.3]{ech-rae} there is a
multiplier $m=(m_C,m_D)$ with $m_C=\psi^{-1}$, and then ${}_C\langle
x,m\rangle=\psi(x)$ for all $x\in X$. Thus for $x\in X$ we have
\begin{align*}
x&=\psi^{-1}\circ\psi(x)=\psi^{-1}\big({}_C\langle x,m\rangle\big)
=m_C\big({}_C\langle x,m\rangle\big)\\
&={}_C\langle x,m\rangle\cdot m
=x\cdot\langle m,m\rangle_D,
\end{align*}
and for $c\in C$ we have
\[
c=\psi\circ\psi^{-1}(c)={}_C\langle \psi^{-1}(c),m\rangle
={}_C\langle c\cdot m,m\rangle=c{}_C\langle m,m\rangle,
\]
from which (\ref{mult}) follows. It is easy to check that, given $m\in
M(X)$ satisfying (\ref{mult}), $\psi:x\mapsto {}_C\langle x,m\rangle$ is
a Hilbert module isomorphism with inverse $c\mapsto c\cdot m$.

For part (ii), define $\alpha = \psi \circ \phi^{-1}$, and note that
$\alpha(d)={}_C\langle m\cdot d,m\rangle$. Then 
$\alpha$ is a linear isomorphism  of
$D$ onto
$C$; to see it  is in fact a \cstar isomorphism requires  only 
calculations using the properties of $m$. For example, for $a,b \in D$ 
we have
\begin{equation*}
\begin{split}
{}_C\langle m\cdot ab,m \rangle
& = {}_C\langle m\cdot a \langle m,m \rangle_D \,b,m \rangle
={}_C\langle{}_C\langle m\cdot a,m \rangle m\cdot b,m \rangle \\
& = {}_C\langle m\cdot a,m \rangle \,{}_C\langle m\cdot b,m \rangle. \qed
\end{split}
\end{equation*}
\noqed
\end{proof}

In the above proof we could have constructed $\alpha$ directly 
from $\psi$: a Hilbert module isomorphism induces an isomorphism of 
the algebras of compact operators, and hence $D=\K_C(X)\cong\K_C(C)
= C$. 

Proposition \ref{module iso} simplifies part of \cite[Section
4]{exe:partial}, where   both $\psi$
and $\alpha$ are postulated, and 
$\psi$ is postulated to be a $C-D$ bimodule isomorphism rather than just a 
left $C$-module isomorphism. Our approach also justifies Exel's feeling
that one of his results 
\cite[Proposition 4.13]{exe:partial} is a kind of ``partial multiplier'' 
property: he is really using a module 
multiplier in the sense of \cite{ech-rae}. 

Now suppose $B$ is a \cstar algebra and $X$ is a closed subspace of $B$ 
such that $XX^*X \subset X$. Let $C = XX^*$ (caution: we use Exel's 
 convention that this denotes the \emph{closed} linear 
span of the set of products!) and $D = X^*X$. Then $C$ and $D$ are \cstar 
subalgebras of $B$, and $X$ is a $C-D$ imprimitivity bimodule. Let $p$ and 
$q$ be the identities of $C^{**}$ and $D^{**}$, respectively, regarded as 
projections in $B^{**}$. \cite[Proposition 2.4]{ech-rae} implies that we 
can identify $M(X)$ with 
\[ \{ m \in pB^{**}q \mid Cm \cup mD \subset X \}, \]
in such a way that the module actions are given by multiplication in
$B^{**}$, and the inner products are given by, for example,
${}_C\langle m,n\rangle=mn^*$. Fortunately, there is no
ambiguity between the various meanings of $m^*$: the adjoint
$m^*$ of $m$ in $B^{**}$ implements the adjoint of the module
homomorphism $c\mapsto c\cdot m$, which is given by $x\mapsto
{}_C\langle x,m\rangle=xm^*$. Finally we observe that, if
$\rho$ is any faithful nondegenerate representation of
$B$ on 
$\H$, then the canonical extension of $\rho$ to a normal representation of 
$B^{**}$ on $\H$ maps $M(X)$ isomorphically onto 
\[ \{ x \in \rho(p)\B(\H)\rho(q) \mid
\rho(C)x \cup x\rho(D) \subset \rho(X) \}. \]

\section{Crossed products and the dual coaction}\label{cps}

Suppose $\alpha$ is a partial action of a group $G$ on a $C^*$-algebra
$A$, and $(\pi,u)$ is a covariant representation of $(A,G,\alpha)$ on a
Hilbert space $\H$. Let
\[ C^*(\pi,u) = \overline{\sum_{s\in G} \pi(D_s)u_s}. \]
A short computation using the relation
$\alpha_s(D_{s\inv}D_t)\subset D_{st}$ shows  that for
$s,t \in G$, $a \in D_s$ and $b \in D_t$ we have 
\[ \pi(a)u_s\pi(b)u_t
= \pi \circ \alpha_s(\alpha_{s\inv}(a)b)u_{st} \]
and
\[ (\pi(a)u_s)^* = \pi \circ \alpha_{s\inv}(a^*)u^*_s, \]
so the closed subspace $C^*(\pi,u)$ of $B(\H)$ is a \cstar algebra,
which we call the \emph{$C^*$-algebra of the covariant representation
$(\pi,u)$}. 

We would like to define the crossed product
$A
\times_\alpha G$ to be the
\cstar  algebra of a universal covariant representation $(\pi,u)$. The
image
$\pi(A)$ is a 
\cstar subalgebra of $C^*(\pi,u)$, but in general (as can be seen from 
\exref{shift} and \propref{multiplier}) the partial isometries $u_s$
need  not be multipliers of  $C^*(\pi,u)$. So to get a suitable
universal covariant representation, we shall have to work in the double
dual
$(A\times_\alpha G)^{**}$, and  we shall have to
construct the algebra
$A
\times_\alpha G$ as an enveloping algebra.

Following McClanahan's development, let $L_c$ denote the vector space of
 functions $f:G\to A$ of finite support such that
$f(s) \in D_s$ for all $s \in G$. For $a \in D_s$, let
\[ F(a,s)(t) =
\begin{cases}
a & \text{if } t = s, \\
0 & \text{otherwise.}
\end{cases} \]
Then $L_c$ is the linear span of the $F(a,s)$. (McClanahan writes $a 
\delta_s$ for our $F(a,s)$; we avoid the notation $\delta_s$ because
it has other connotations for coaction freaks, and call them 
$m_s$ instead. The precise meaning of $m_s$ will be made clear shortly.)
 The \star algebra
structure of $L_c$ is determined by 
\begin{align*}
F(a,s)F(b,t) &= F(\alpha_s(\alpha_s\inv(a)b),st), \ \text {and}\\
F(a,s)^* &= F(\alpha_s\inv(a^*),s\inv).
\end{align*}
Note that $A$ embeds as a \star subalgebra of $L_c$ via $a \mapsto F(a,e)$, 
so any \star representation $\Pi$ of $L_c$ on a Hilbert space restricts to 
a representation of the \cstar algebra $A$. Hence 
\begin{equation}\label{bounded}
\begin{split}
\| \Pi(F(a,s)) \|^2
& = \| \Pi(F(a,s))^*\Pi(F(a,s)) \| \\
& = \| \Pi(F(\alpha_s\inv(a^*),s\inv)F(a,s)) \| \\
& = \| \Pi(F(\alpha_s\inv(a^* a),e)) \|
\le \| \alpha_s\inv(a^*a) \|
= \| a \|^2.
\end{split}
\end{equation}
We deduce that the greatest \cstar seminorm on $L_c$ is finite (it is in 
fact a norm, although we do not need this), so it makes sense to define the 
crossed product $A \times_\alpha G$ as the \cstar completion of $L_c$.
The double dual
$A^{**}$ embeds naturally in $(A 
\times_\alpha G)^{**}$, so  the projections $p_s$ are naturally
identified  with (no longer central) projections in $(A \times_\alpha
G)^{**}$. 

Let
\[ X_s = \{ F(a,s) \mid a \in D_s \}. \]
The calculation (\ref{bounded}) shows that this is a closed subspace of $A
\times_\alpha G$, and 
$ X_sX^*_s = D_s$ and $X^*_sX_s = D_{s\inv}$ as subalgebras of $A\subset
A\times_\alpha G$, so $X_s$ is a $D_s - D_{s\inv}$ imprimitivity
bimodule. The discussion at the end of \S2 shows that 
\begin{equation}\label{multxs}
 M(X_s) = \{ m \in p_s(A \times_\alpha
G)^{**}p_{s\inv}
\mid D_sm \cup mD_{s\inv} \subset X_s \}.
\end{equation}

A calculation shows that the maps $l:D_s\to X_s$, $r:D_{s^{-1}}\to X_s$
defined by
\[
l(a)= F(a,s),\ \ 
r(b)= F(\alpha_s(b),s)
\]
satisfy $l(a)\cdot b=a\cdot r(b)$, and hence define a multiplier $m_s:=
(l,r)$ of
$M(X_s)$ \cite[Lemma 1.4]{ech-rae}. The identification (\ref{multxs})
allows us to view
$m_s$ as an element of 
$p_s(A \times_\alpha G)^{**}p_{s^{-1}}$, which by definition satisfies
\[
a\cdot m_s =F(a,s)\  \text{for  $a \in D_s$, and }
m_s\cdot b=F(\alpha_s(b),s)\ \text{for $b \in D_{s\inv}$}.
\]
Recall that the $M(D_s)$-valued inner product is given on $p_s(A \times_\alpha
G)^{**}p_{s\inv}$ by ${}_{D_s}\langle m,n\rangle=mn^*$; thus for $a\in
D_s$ we have
\[
a(m_sm_s^*)=a{}_{D_s}\langle m_s,m_s\rangle=
{}_{D_s}\langle a\cdot m_s,m_s\rangle=
{}_{D_s}\langle F(a,s),m_s\rangle.\]
It follows from (\ref{ipsonm(x)}) that ${}_{D_s}\langle F(a,s),m_s\rangle$
is given in terms of the adjoint of $l:D_s\to X_s$ by $l^*(F(a,s))$,
which by an easy calculation is seen to be $a$. Thus $a(m_sm_s^*)=a$ for
all $a\in D_s$, and $m_sm_s^*=1_{M(D_s)}=p_s$. Similiarly, we can verify
that $m_s^*m_s=p_{s^{-1}}$. The multipliers $m_s$ therfore have the
property (\ref{mult}) of Lemma \ref{module iso}, and induce 
isomorphisms $\ad m_s$ of $D_{s^{-1}}$ onto $D_s$, characterised by $\ad
m_s(a)={}_{D_s}\langle m_s\cdot a,m_s\rangle$. Another calculation shows
that 
\[
\ad m_s(a)={}_{D_s}\langle m_s\cdot a,m_s\rangle=l^*(m_s\cdot a)=
l^*\big(F(\alpha_s(a),s)\big)=\alpha_s(a),
\]
so the elements $m_s$ of $(A\times_\alpha G)^{**}$ are partial
isometries implementing the partial automorphisms $\alpha_s$.

We claim that the inclusion $\iota \colon A \hookrightarrow A 
\times_\alpha G$ and $m$ form a covariant representation $(\iota,m)$ of 
$(A,G,\alpha)$ in $(A \times_\alpha G)^{**}$. To see this, let  $a \in
D_sD_{st}$, and let $e_i$ be an approximate identity  for $D_t$. Then in
the weak* topology of $(A \times_\alpha G)^{**}$, we have
\begin{equation*}
\begin{split}
am_sm_t
& = \lim am_se_im_t
= \lim \alpha_s(\alpha_{s\inv}(a)e_i)m_{st} \\
& = \alpha_s(\alpha_{s\inv}(a))m_{st}
= am_{st}.
\end{split}
\end{equation*}
It now follows from \lemref{altcovrep} that  $(\iota,m)$ is covariant,
as claimed.

Next, let $(\pi,u)$ be a
covariant representation of $(A,G,\alpha)$.  Then $(\pi,u)$ defines a
representation of $L_c$, which extends to a representation
$\pi\times u$ of $A\times_\alpha G$ by definition of the enveloping norm,
and hence also to a normal representation, also denoted $\pi\times u$, 
of $(A
\times_\alpha G)^{**}$. Applying $\pi\times u$ to $\iota(a)=F(a,e)$ gives
$\pi(a)$, and for
$a\in D_{s^{-1}}$ we have
\begin{equation*}
\begin{split}
\pi\times u(m_s)(\pi(a)h)
& =\pi\times u\big((m_sa)h\big)=\pi\times u\big(F(\alpha_s(a),s)\big)h \\
&
=\pi(\alpha_s(a))u_sh=u_s\pi(a)u_s^*u_sh=u_s\pi(a)\pi(p_s)h\\
&=u_s\pi(ap_s)h=u_s(\pi(a)h),
\end{split}
\end{equation*}
so $\pi\times u(m_s)=u_s$. Thus
\begin{equation}
\label{separatecovariant}
(\pi \times u) \circ \iota = \pi
\spacetext{and} (\pi \times u) \circ m = u. 
\end{equation}
If we want to avoid extending to $(A \times_\alpha 
G)^{**}$, we need to rewrite \eqref{separatecovariant} as
\[ (\pi \times u)(am_s) = \pi(a)u_s
\spacetext{for} a \in D_s. \]  
That covariant representations extend to $A\times_\alpha G$ in this way
characterises the crossed product:

\begin{prop}
Let $(\pi,u)$ be a covariant representation of $(A,G,\linebreak[0]\alpha)$.
Then $\pi \times u$ is faithful if and only if for every covariant 
representation $(\rho,v)$ there exists a homomorphism $\theta \colon 
C^*(\pi,u) \to C^*(\rho,v)$ such that
\[ \theta(\pi(a)u_s) = \rho(a)v_s \spacetext{for} a \in D_s. \]
\end{prop}

\begin{proof}
Applying the hypothesis to the canonical covariant representation
$(\iota,m)$  of $(A,G,\alpha)$ in $(A \times_\alpha G)^{**}$ gives an
inverse $\theta$ for
$\pi \times u$. 
\end{proof}

Although this is not a deep result, it 
does give the easiest way of recognising the full crossed product.
However,  it is difficult to see how to turn this into a convenient
categorical definition of  the crossed product, as is often done for
ordinary actions. For  example, even if $\pi \times u$ is faithful on $A
\times_\alpha G$, we  do not know if the canonical extension to a normal
representation of $(A 
\times_\alpha G)^{**}$ is faithful on the \cstar algebra generated by 
$A$ and $m_G$. Thus for all we know, $C^*(A \cup m_G)$ and 
$C^*(\pi(A) \cup u_G)$ could be essentially different even though $A 
\times_\alpha G$ and $C^*(\pi,u)$ are isomorphic. For a partial solution 
of this conundrum, see \propref{faithfulu} below. 

\begin{prop}
Let $\alpha$ be a partial action of a discrete group $G$ on a
$C^*$-algebra $A$. Then there is a unique  coaction $\hat\alpha$ of $G$ on
$A
\times_\alpha G$ such that
$\hat\alpha(am_s)  = am_s \otimes s$ for $a \in D_s$. 
\end{prop}

\begin{proof}
Define maps $\pi$ and $u$ from $A$ and $G$, respectively, to $(A 
\times_\alpha G)^{**} \otimes C^*(G)$ by 
\[ \pi(a) = a \otimes 1 \spacetext{and} u_s = m_s \otimes s. \]
Then $u$ is a partial representation (being a tensor product of two such). 
We have 
\begin{gather*}
u_su^*_s
= m_sm^*_s \otimes 1 = p_s \otimes 1 = \pi(p_s), \\
u_su_tu^*_tu^*_s
= m_sm_tm^*_tm^*_s \otimes stt\inv s\inv
= m_sm_tm^*_{st} \otimes 1
= u_su_tu^*_{st},
\end{gather*}
and
\begin{equation*}
\begin{split}
\ad u_s \circ \pi(a)
& = \ad m_s(a) \otimes 1 \\
& = \alpha_s(a) \otimes 1
= \pi \circ \alpha_s(a) \righttext{for} a \in D_{s\inv},
\end{split}
\end{equation*}
so $(\pi,u)$ is a covariant representation of $(A,G,\alpha)$ by Lemma
\ref{altcovrep}(ii). Clearly
$C^*(\pi,u) \subset (A \times_\alpha G) \otimes C^*(G)$.

Define $\theta \colon (A \times_\alpha G) \otimes C^*(G) \to A 
\times_\alpha G$ by $\theta(x \otimes s) = x$ for $s \in G$. This is
well-defined on the  minimal tensor product because it is the 
tensor product of the identity homomorphism of $A \times_\alpha G$ and 
the augmentation  representation of $C^*(G)$. We have 
\[ \theta \circ (\pi \times u)(am_s) = \theta(am_s \otimes s) = am_s
\spacetext{for} a \in D_s. \]
Hence, $\hat\alpha = \pi \times u$ is a faithful homomorphism of $A 
\times_\alpha G$. Nondegeneracy  of $\hat\alpha$ as a homomorphism into
$(A 
\times_\alpha G) \otimes C^*(G)$ and the coaction identity are obvious. 
\end{proof}

As a first application, we use the dual coaction $\hat\alpha$ to show
that in a faithful representation $\pi\times u$, $u$ is as faithful as
$m$ is: 

\begin{prop}
\label{faithfulu}
If $\pi \times u$ is a faithful representation of $A \times_\alpha G$, 
then  $u_s=u_t$ implies $m_s=m_t$.
\end{prop}

\begin{proof}
Suppose $u_s = u_t$, and suppose first that $D_sD_t\not=\{0\}$. Then for
any nonzero
$a
\in D_sD_t$ we have
\[ (\pi \times u)(am_s) = \pi(a)u_s = \pi(a)u_t
= (\pi \times u)(am_t), \]
so $am_s = am_t$ by hypothesis. Since $a\not=0$ implies
$am_s\not=0$, applying the dual coaction 
$\hat\alpha$ gives $am_s \otimes s = am_t \otimes t$, so $s$ and $t$ 
are linearly dependent elements of $C^*(G)$, forcing $s=t$ and $m_s=m_t$.
If
$D_s\cap D_t=D_sD_t=\{0\}$, then $p_sp_t=0$, and 
\[u_su_s^*=u_su_s^*u_tu_t^*=\pi(p_s)\pi(p_t)=\pi(p_sp_t)=0,
\]
forcing $u_s=0$. But then $\pi\times u(am_s)=\pi(a)u_s=0$ for all $a\in
D_s$, and because $\pi\times u$ is faithful this implies $D_s=\{0\}$ and
$m_s=0$. Simlarly, $u_t=u_s=0$ implies $m_t=0$, so again we
have $m_s=m_t$, as required.
\end{proof}

Recall from \cite{qui:disc} that if $\delta$ is a coaction of $G$ on a
\cstar  algebra $B$, then for $s\in G$ the associated spectral subspace
is 
\[ B_s  =\{ b \in B \mid \delta(b) = b \otimes s \}.\]
If $\chi_s$ is the characteristic function of $\{s\}$, regarded as an 
element of $B(G) = C^*(G)^*$, then $\delta_s = (\iota \otimes \chi_s)
\circ \delta$ is a projection of $B$ onto $B_s$.

\begin{prop}
The spectral subspaces for the dual coaction $\hat\alpha$ on a partial
crossed product
$A\times_\alpha G$ are given by
$(A
\times_\alpha  G)_s = D_sm_s = m_sD_{s\inv}$. 
\end{prop}

\begin{proof}
Let $\hat\alpha_s \colon A \times_\alpha G \to (A \times_\alpha G)_s$ be 
the canonical projection. Clearly $D_sm_s \subset (A \times_\alpha 
G)_s$. On the other hand, any  $x \in (A \times_\alpha G)_s$ can be 
approximated by a finitely nonzero sum $\sum_t a_tm_t$, and then 
\[ x = \hat\alpha_s(x) \approx \hat\alpha_s(\sum_t a_tm_t) = a_sm_s. \]
Hence $(A \times_\alpha G)_s \subset D_sm_s$, proving the first 
equality. The second equality follows from the covariance of $(\iota,m)$:
\[ m_sD_{s\inv} = m_sD_{s\inv}p_{s\inv} = m_sD_{s\inv}m^*_sm_s
= \alpha_s(D_{s\inv})m_s = D_sm_s. \qed \]
\noqed
\end{proof}

\begin{prop}
\label{multiplier}
For $s \in G$, consider the following conditions\textup:
\begin{enumerate}
\item $p_s \in M(A)$\textup;
\item $D_s = Ap_s$\textup;
\item $M(D_s) \subset M(A)$\textup;
\item $m_s \in M(A \times_\alpha G)$.
\end{enumerate}
Conditions \textup{(i)--(iii)} are equivalent, and are implied by 
\textup{(iv)}. Moreover, if \textup{(i)} holds for both $s$ and $s\inv$, 
then \textup{(iv)} holds as well. 
\end{prop}

\begin{proof}
That (i) implies (ii) must be a well-known general 
fact about ideals of \cstar algebras, but we lack a reference. Let
$A$ act via its universal representation  on $\H$. It
suffices to show that any state
$\omega$ of $A$  annihilating $D_s$ also annihilates $Ap_s$. There exists
$\xi \in \H$  such that $\omega(a) = (a\xi,\xi)$. Since $p_s$ is in the
weak* closure  of $D_s$, $(p_s\xi,\xi)=0$. This forces $p_s\xi = 0$, so
for any $a \in A$ we have $\omega(ap_s) = (ap_s\xi,\xi) = 0$, as
required.

The chain (ii) implies (iii) implies (i) is routine. 

Assuming (iv), we have $p_s = m_sm^*_s \in M(A \times_\alpha G)$ 
also. Since $\hat\alpha(p_s) = p_s \otimes 1$, we get 
\[ p_s \in M(A \times_\alpha G)_e  = M((A \times_\alpha G)_e)
= M(D_em_e) = M(A). \]

Finally, assume (i) holds for both $s$ and $s\inv$. Since $A 
\times_\alpha G = \overline{\sum_t D_tm_t}$ and $m^*_s = m_{s\inv}$,
(iv)  follows from the following computation for $a\in D_t$: 
\begin{alignat*}{2}
am_tm_s
& = m_t\alpha_{t\inv}(a)p_sm_s \\
& = \alpha_t(\alpha_{t\inv}(a)p_s)m_tm_s
&& \righttext{since} \alpha_{t\inv}(a)p_s \in D_{t\inv}D_s \\
& = \alpha_t(\alpha_{t\inv}(a)p_s)m_{ts}
&& \righttext{since} \alpha_t(\alpha_{t\inv}(a)p_s) \in D_tD_{ts} \\
& \in D_{ts}m_{ts} \subset A \times_\alpha G. \qed
\end{alignat*}
\noqed
\end{proof}

McClanahan \cite{mcc} constructs a regular covariant representation
$(\pi^r,u^r)$, of $(A,G,\alpha)$. We give a  description which is more
convenient for our purposes. For $s \in G$ let 
$\bar\alpha_s \colon A \to M(D_s)$ be the canonical homomorphism
extending 
$\alpha_s \colon D_{s\inv} \to D_s$; also let $\chi_s$ be the 
characteristic function of $\{s\}$, and  $\lambda$ be the left 
regular representation of $G$. Let $\pi$ be a faithful and nondegenerate
representation of
$A$ on $\H$. Then
$(\pi^r,u^r)$ acts on
$\H
\otimes l^2(G)$, and is determined by 
\[ \pi^r \times u^r(am_s)
= \sum_t \bar\alpha_{t\inv}(a) \otimes \chi_t\lambda_s
\spacetext{for} a \in D_s, \]
where the sum converges in the strong* topology. The \emph{reduced 
crossed product} of $(A,G,\alpha)$ is 
$A \times_{\alpha,r} G = C^*(\pi^r,u^r)$. 

The following result characterises the reduced crossed product as the 
canonical image of $A \times_\alpha G$ in the multipliers of the double 
crossed product $(A \times_\alpha G) \times_{\hat\alpha} G$: 

\begin{prop}
Let $\alpha$ be a partial action of a discrete group $G$ on a
$C^*$-algebra $A$, and let $j_{A\times_\alpha G}$ be the canonical
embedding of $A\times_\alpha G$ in the crossed product by the dual
coaction. Then there is an isomorphism
$\theta
\colon j_{A
\times_\alpha G}(A 
\times_\alpha G) \to A \times_{\alpha,r} G$ such that 
\[ \theta \circ j_{A \times_\alpha G} = \pi^r \times u^r. \]
\end{prop}

\begin{proof}
Let $M$ be the  representation of $c_0(G)$ by multiplication operators
on 
$l^2(G)$. By \cite[Lemma 2.2]{qui:disc}, the following calculation
shows 
$(\pi^r \times u^r,1 \otimes M)$ is a covariant representation of $(A 
\times_\alpha G,G,\hat\alpha)$: for $a \in D_s$ 
\begin{equation*}
\begin{split}
(\pi^r \times u^r)(am_s)(1 \otimes \chi_t)
& = \sum_r \bar\alpha_{r\inv}(a) \otimes \chi_r\lambda_s\chi_t \\
& = \sum_r \bar\alpha_{r\inv}(a) \otimes\chi_r\chi_{st}\lambda_s \\
& = (1 \otimes \chi_{st})(\pi^r \times u^r)(am_s).
\end{split}
\end{equation*}
Since $(\pi^r \times u^r) |_{(A \times_\alpha G)^{\hat\alpha}} = (\pi^r 
\times u^r) |_A = \pi^r$ is faithful, \cite[Proposition 2.18]{qui:disc} 
shows $\ker(\pi^r \times u^r) = \ker j_{A \times_\alpha G}$, and the 
result follows. 
\end{proof}

\begin{rem}
The above proposition shows $A \times_{\alpha,r} G$ is independent up to 
isomorphism of the choice of faithful representation of $A$, so
\cite[Proposition 3.4]{mcc} is a corollary. By \cite[Proposition  2.8
(i)]{qui:full} and
\cite[Proposition 2.6 and Theorem 4.1  (2)]{rae:full} $\ker j_{A
\times_\alpha G} = \ker(\iota \otimes \lambda) 
\circ \hat\alpha$. Thus we also obtain alternative proofs of 
\cite[Lemma 4.1]{mcc} and the half of \cite[Proposition 4.2]{mcc} 
stating that if $G$ is amenable then $A \times_\alpha G = A 
\times_{\alpha,r} G$. 
\end{rem}

A coaction $\delta$ of $G$ on a \cstar algebra $B$ is called
\emph{normal} \cite{qui:fullred} if $j_B$ (or equivalently, $(\iota
\otimes 
\lambda) \circ \delta$) is faithful. The coaction $\ad(j_G 
\otimes \iota)(w_G)$ on $j_B(B)$ is always normal, and has the same
covariant  representations and crossed product as $\delta$
\cite[Proposition  2.6]{qui:fullred}; it is called the
\emph{normalization} of $\delta$,  and denoted $\delta^n$. The previous
proposition allows us to view the  normalization $\hat\alpha^n$ of
the dual coaction $\hat\alpha$ as a  coaction on the reduced crossed
product $A
\times_{\alpha,r} G$. 

\begin{cor}
\label{regcovrep}
Let $(\pi,u)$ be a covariant representation of $(A,G,\alpha)$. Then 
$\ker(\pi\times u) = \ker(\pi^r \times u^r)$ if and only if $\pi$ is 
faithful and there is a normal coaction $\delta$ of $G$ on $C^*(\pi,u)$ 
with $\delta \circ (\pi \times u) = ((\pi \times u) \otimes \iota) \circ 
\hat\alpha$. 
\end{cor}

\begin{proof}
This is immediate from the Proposition and \cite[Corollary 
2.19]{qui:disc}. 
\end{proof}

\section{Landstad duality}

 Let
$\delta$ be a coaction of the discrete group $G$ on a \cstar algebra
$B$. For $s\in G$ let $B_s:=\{b\in
B\mid\delta(b)=b\otimes s\}$ be the spectral subspace, and let
$D_s=B_sB^*_s:=\overline{\text{sp}}\{bc^*:b,c\in B_s\}$. Then $D_s$ is an
ideal of
$D_e=B_e$, and
$D_{s\inv}=B^*_sB_s$. Let $p_s$ denote the identity of $D^{**}_s$
regarded as a projection in $B^{**}$.
$B_s$ is a $D_s-D_{s\inv}$ imprimitivity bimodule with inner products
$_{D_s}\langle x,y\rangle=xy^*$ and $\langle
x,y\rangle_{D_{s\inv}}=x^*y$.
By the discussion at the end of \secref{mbimods}, the multiplier bimodule
can be identified as
\[ M(B_s) = \{ b \in p_sB^{**}p_{s\inv}
\mid D_sb \cup bD_{s\inv} \subset B_s \}. \]
Fortunately, when $s=e$ this coincides with the usual multiplier algebra 
$M(B_e)$ of $B_e$. 

The following result is  Landstad duality for partial actions. 
Condition \eqref{module iso condition} below was motivated by 
\cite[Proposition 4.16]{exe:partial}. 

\begin{thm}
\label{duality}
Let $\delta$ be a normal coaction of a discrete group $G$ on a \cstar
algebra $B$.  The following are equivalent\textup:
\begin{enumerate}
\item there is a partial action $\alpha$ of $G$ on a \cstar algebra $A$ 
such that $(B,\delta)$ is isomorphic to 
$(A\times_{\alpha,r}G,\hat\alpha^n)$\textup; 
\item there is a partial representation $m$ of $G$ in $B^{**}$ such
that
\begin{equation}
\label{multiplier condition}
m_s\in M(B_s)\spacetext{and}m_sm^*_s=p_s\spacetext{for}s\in G\textup;
\end{equation}
\item there is a collection $\psi_s\colon B_s\to D_s$ 
of left Hil\-bert $D_s$-mod\-ule isomorphisms
such that
\begin{equation}
\label{module iso condition}
\psi_{st}(xy)=\psi_s(x\psi_t(y))\spacetext{for}x\in B_s,y\in B_t.
\end{equation}
\end{enumerate}
\end{thm}

\begin{proof}
The construction in \secref{cps} shows that (i) implies (ii). We next
show that (ii) implies (i). Since $m$ is a partial representation, we have
$m^*_sm_s=p_{s\inv}$. By
\lemref{module iso}, there are isomorphisms
$\alpha_s:=\ad m_s\colon D_{s\inv}\to D_s$; we claim that $\alpha$ is a
partial action of $G$ on $B_e$. Clearly $D_e=B_e$. We must show
\[ \alpha_s(D_{s\inv}D_t)\subset D_{st} \]
and
\[ \alpha_s\alpha_t=\alpha_{st}\spacetext{on}D_{t\inv}D_{t\inv s\inv}. \]
For the first,
\[ \alpha_s(D_{s\inv}D_t)
=m_sD_{s\inv}D_tD_{s\inv}m^*_s
\subset B_sB_tB^*_tB^*_s
\subset B_{st}B^*_{st}
=D_{st}. \]
For the second, since $m$ is a partial representation,
\eqref{multiplier condition} and \lemref{altpartrep} imply
$am_sm_t=am_{st}$ for all $a\in D_sD_{st}$, or equivalently
$m_{st}a=m_sm_ta$ for all $a\in D_{t\inv}D_{t\inv s\inv}$. Thus for
such $a$ we have
\[ \alpha_{st}(a)=m_{st}am^*_{st}=m_sm_tam^*_tm^*_s=\alpha_s\alpha_t(a),
\]
as claimed.

The pair $(\iota,m)$
is a covariant representation of $(B_e,G,\alpha)$, and we
have $C^*(\iota,m)=B$ because $B=\overline{\sum_sB_s}$. For $a\in
D_s$ we have
\begin{equation*}
\begin{split}
\delta\circ(\iota\times m)(am_s)
& =\delta(am_s)
=am_s\otimes s
=(\iota\times m)\otimes\iota(am_s\otimes s)\\
& =((\iota\times m)\otimes\iota)\circ\hat\alpha(am_s);
\end{split}
\end{equation*}
since $\iota$ is faithful, (i) now follows from Corollary \ref{regcovrep}

Now we show that (ii) implies (iii). By \lemref{module iso},
$\psi_s(x)=xm^*_s$ defines a left Hilbert module isomorphism
$\psi_s\colon B_s\to D_s$. Let $x\in B_s$ and $y\in B_t$. Then there
exist $a\in D_{s\inv}$ and $b\in D_t$ with $x=m_sa$ and $y=bm_t$. We
compute:
\begin{equation*}
\begin{split}
\psi_{st}(xy)
& =m_sabm_tm^*_{st}
=m_sabm^*_sm_sm_tm^*_{st}\righttext{since}ab\in D_{s\inv}\\
& =m_sabm^*_sm_sm_tm^*_tm^*_s
=m_sabm_tm^*_tm^*_s\\
& =xym^*_tm^*_s
=\psi_s(x\psi_t(y)),
\end{split}
\end{equation*}
giving \eqref{module iso condition}.

Finally, to see that (iii) implies (ii), \lemref{module iso} gives $m_s\in
M(B_s)$ such that
\[ m_sm^*_s=p_s\spacetext{and}m^*_sm_s=p_{s\inv} \]
and it remains to verify \eqref{hom}.  But \eqref{module iso
condition} implies
\[ am_sbm_tm^*_{st}=am_sbm_tm^*_tm^*_s\spacetext{for}a\in D_s,b\in D_t, \]
and letting $a$ and $b$ run separately through bounded approximate
identities for $D_s$ and $D_t$ gives
$m_sm_tm^*_{st}=m_sm_tm^*_tm^*_s$.
\end{proof}

Landstad \cite[Theorem 3]{lan:dual} originally characterised reduced
crossed products by ordinary actions of a locally compact group. When
the group is discrete, Landstad's characterisation is the special case
of the above theorem in which $p_s=1$ for all $s\in G$: $\delta$
is equivariantly isomorphic to the dual coaction on a reduced crossed
product by an ordinary action of $G$ if and only if there is a
homomorphism $s\mapsto m_s$ of $G$ into $UM(B)$
such that $\delta(m_s)=m_s\otimes s$ for all $s\in G$.

\section{Partial actions of $\F_n$}

McClanahan \cite[Example 2.4]{mcc} and (for the case $n=1$) Exel
\cite{exe:partial} show that certain partial actions of the free group $\F_n$
can be reconstructed from the generators. We shall show that, for such
partial actions of $\F_n$, our Landstad duality (\thmref{duality}) can be
recast in terms of the spectral subspaces of the generators. This will
give both a generalisation and a simplification of Exel's
characterisation \cite[Theorem 4.21]{exe:partial} of crossed products by
certain partial actions of $\Z=\F_1$.

Throughout these last two sections we shall denote by
$g_1,\dots,g_n$ a fixed set of generators for the free group $\F_n$. A
word in $\F_n$ is \emph{reduced} if it is the identity or a product
$s_1s_2\cdots s_k$ in which each $s_i$ is either $g_j$ or $g_j^{-1}$ for
some $j$, and no cancellation is possible. When we say $s_1s_2\cdots
s_k$ is a reduced word it is implicit that each $s_i$ has the form
$g_j^{\pm 1}$.

\begin{defn}
A partial action $\alpha$ of $\F_n$ is \emph{multiplicative} if
 for every reduced word 
$s_1\cdots s_k$, we have
$\alpha_{s_1\cdots s_k}=\alpha_{s_1}\cdots\alpha_{s_k}$.
\end{defn}

Thus for any multiplicative partial action we have
$\alpha_{st}=\alpha_s\alpha_t$ whenever
$s,t$ are words for which there is no cancellation possible in $st$. The
key issue is that the domains of definition must coincide, and  we shall
see in the next Lemma that it is relatively easy to decide whether this
happens.

\begin{lem}
\label{elementaryprop}
A partial action $\alpha$ of $\F_n$ is multiplicative if and only if
$D_{s_1s_2\cdots s_k}\subset D_{s_1}$ for every reduced word $s_1\cdots
s_k$.
\end{lem}

\begin{proof}
The ``only if'' direction is clear. For the other direction, it suffices
to show that $\alpha_{s_1\cdots s_k}=\alpha_{s_1}\alpha_{s_2\cdots
s_k}$. Write $s=s_1$, $t=s_2\cdots s_k$. Then because $\alpha_{st}$ is
an injective map which extends $\alpha_s\alpha_t$, it is enough to show
they have the same range. But
\begin{equation*}
\begin{split}
\operatorname{range}\alpha_{st}
&=D_{st}=D_sD_{st}\ \mbox{ by hypothesis}\\
&=\alpha_s(D_s^{-1}D_t)\ \mbox{ by Lemma \ref{domains}}\\
&=\alpha_s\alpha_t(D_{t^{-1}s^{-1}}D_{t^{-1}})
\ \mbox{ by Lemma \ref{domains}}\\
&=\operatorname{range}\alpha_s\alpha_t,
\end{split}
\end{equation*}
as required.
\end{proof}

\begin{lem}
If $\alpha$ is a multiplicative partial action of $\F_n$ on a
$C^*$-algebra $A$, then the partial representation $m$ of $\F_n$ in
$(A\times_\alpha\F_n)^{**}$ satisfies $m_{s_1\cdots s_k}=m_{s_1}\cdots
m_{s_k}$ for every reduced word $s_1\cdots s_k$.
\end{lem}

\begin{proof}
Again write $s=s_1$, $t=s_2\cdots s_k$, and it is enough to show that
the partial isometries $m_{st}$ and $m_sm_t$ have the same range
projection. But since $D_{st}=D_sD_{st}=\alpha_s(D_{s^{-1}}D_t)$, we
have 
\[
p_{st}=\alpha_s(p_{s^{-1}}p_t)=m_sp_{s^{-1}}p_tm_s^*
=m_sm_tm_t^*m_s^*.\qed
\]
\noqed
\end{proof}

In the above situation, $D_{s_1\cdots s_k}
=\alpha_{s_1}(D_{s\inv_1}D_{s_2\cdots s_k})$, so $\alpha$ is a free
product of $n$ partial actions of $\Z$ in the sense of McClanahan
\cite[Example 2.3]{mcc}. In fact, a partial action of $\F_n$ is
multiplicative if and only if it is the free product of $n$
multiplicative partial actions of $\Z$.   This brings up a minor
inconsistency between the partial actions of Exel and McClanahan's
partial actions of
$\Z$: a partial action of $\Z$ in McClanahan's sense \cite{mcc} is a
partial action in Exel's sense \cite{exe:partial} if and only if it is
multiplicative. Examples of nonmultiplicative partial actions of $\Z$ are
easy to come by:

\begin{ex}
Suppose $\beta$ is an (ordinary) action of
$\Z$ on $A$. Define ideals $\{D_n\}_{n\in\Z}$ of $A$ by
\[ D_n=
\begin{cases}
A&\text{if $n$ is even},\\
\{0\}&\text{if $n$ is odd}.
\end{cases} \]
Then define $\alpha_n=\beta_n|D_{-n}$. To see that $\alpha$ is a partial
action, the only nontrivial condition is
$\alpha_n(D_{-n}D_k)\subset D_{n+k}$. This is trivially satisfied when
$n+k$ is even, and if $n+k$ is odd then $n$ or $k$ is odd, and
$\alpha_n(D_{-n}D_k)=\{0\}$. This partial action is not multiplicative
since
$D_2=A\not\subset\{0\}=D_1$.
\end{ex}

It is easy to check whether a partial action of $\Z$ is multiplicative.

\begin{lem}
A partial action $\alpha$ of $\Z$ is multiplicative
if and only if $D_n\subset D_1$ for all $n>0$.
\end{lem}

\begin{proof}
It suffices to show that if $D_n\subset D_1$ for all $n>0$
then $D_{-n}\subset D_{-1}$ for all $n>0$. This is proven inductively by
the following computation:
\begin{equation*}
\begin{split}
D_{-n}
& = \alpha_{-n}(D_n)
= \alpha_{-n}(D_1D_n) \\
& = \alpha_{-n}\alpha_1(D_{-1}D_{n-1})
= \alpha_{1-n}(D_{-1}D_{n-1}) \\
& = \alpha_{1-n}(D_{n-1}D_{-1})
= D_{1-n}D_{-n} \subset D_{1-n}.\qed
\end{split}
\end{equation*}
\noqed
\end{proof}

\begin{thm}
\label{elementaryduality}
Let $\delta$ be a normal coaction of $\F_n$ on a \cstar algebra
$B$. Then there is an multiplicative partial action $\alpha$ of $\F_n$
on a \cstar algebra $A$ such that $(B,\delta)$ is isomorphic to
$(A\times_{\alpha,r}\F_n,\hat\alpha^n)$ if and only if
\begin{enumerate}
\item $B$ is generated by $B_e\cup B_{g_1}\cup\cdots\cup B_{g_n}$\textup;
\item for each $s=g_1,\ldots,g_n$ there exists $m_s\in M(B_s)$ such
that
\[m_sm^*_s=p_s \spacetext{and} m^*_sm_s=p_{s\inv}.\]
\end{enumerate}
Moreover, \textup{(ii)} can be replaced by
\begin{enumerate}
\setcounter{enumi}{2}
\item for each $s=g_1,\ldots,g_n$ there is a left Hilbert
$D_s$-module isomorphism $\psi_s\colon B_s\to D_s$.
\end{enumerate}
\end{thm}

\begin{proof}
First of all, \lemref{module iso} tells
us that (ii) is
equivalent to (iii). If $\alpha$ is a multiplicative partial action of
$\F_n$ on
$A$ and $B=A\times_{\alpha,r}\F_n$, we know (ii) holds. To see (i)
it suffices to show that if $s_1\cdots s_k$ is a reduced word then
\begin{equation}
\label{reducedword}
B_{s_1\cdots s_k}=B_{s_1}\cdots B_{s_k},
\end{equation}
and by induction it suffices to show
\[ B_{s_1\cdots s_k}=B_{s_1}B_{s_2\cdots s_k}. \]
Letting $s=s_1$, $t=s_2\cdots s_k$, and using Lemma \ref{elementaryprop},
we have
\begin{equation*}
\begin{split}
B_{st}
& =D_{st}m_{st}
=D_{s}D_{st}m_{st}
=\alpha_{s}(D_{s\inv}D_t)m_{st}\\
& =\alpha_{s}(\alpha_{s\inv}(D_{s})D_t)m_{st}
=D_{s}m_{s}D_tm_t
=B_{s}B_t,
\end{split}
\end{equation*}
as desired.

Conversely, assume that $\delta$ is a normal coaction of $\F_n$ on $B$
satisfying (i) and (ii). We first show that if $s_1\cdots s_k$ is a
reduced word then
\eqref{reducedword} holds again. Every $x\in B_{s_1\cdots s_k}$
is approximated by a sum of terms of the form $x_1\cdots x_j$, with
$x_i\in B_{t_i}$ and $t_i\in\{e,g^{\pm 1}_1,\ldots,g^{\pm 1}_n\}$.
By applying the canonical projection of $B$ onto $B_{s_1\cdots s_k}$ to
this sum, we see that we can assume that each $x_1\cdots
x_j\in B_{s_1\cdots s_k}$. This forces $t_1\cdots t_j=s_1\cdots s_k$.
Since the latter product is reduced, we can insert parentheses in the
former product so that the $i$th chunk of $t$'s multiplies out to
$s_i$. Hence,
\[ x_1\cdots x_j\in B_{s_1}\cdots B_{s_k}, \]
verifying \eqref{reducedword}.

To apply \thmref{duality}, we need a partial representation
$m$ of $\F_n$ in $B^{**}$ such that for each $s\in\F$
\begin{equation}
\label{multiplier condition for s}
m_s\in M(B_s)\spacetext{and}m_sm^*_s=p_s.
\end{equation}
We are given $m_s$ satisfying \eqref{multiplier condition for s}
for $s=g_1,\ldots,g_n$. Define
\begin{gather*}
m_e=1,\qquad m_{g\inv_i}=m^*_{g_i},\qquad\text{and}\\
m_{s_1\cdots s_k}=m_{s_1}\cdots m_{s_k}
\spacetext{for any reduced word}s_1\cdots s_k.
\end{gather*}
For $s=e,g\inv_1,\ldots,g\inv_n$ it is clear that \eqref{multiplier
condition for s} holds.  Note also that $m^*_s=m_{s\inv}$ for all
$s\in\F_n$. To see that \eqref{multiplier condition for s} remains true
for all $s\in\F_n$, by induction it suffices to show that if
it holds for $s=r,t$ and there is no cancellation in $rt$, then
\begin{gather}
m_{rt}\in M(B_{rt});\label{multiplierproduct1}\\
m_{rt}m^*_{rt}=p_{rt}.\label{multiplierproduct2}
\end{gather}
For \eqref{multiplierproduct1} we need
\begin{gather}
\label{multiplierproduct11}
m_{rt}\in p_{rt}B^{**}p_{t\inv r\inv};\\
\label{multiplierproduct12}
D_{rt}m_{rt}\cup m_{rt}D_{t\inv r\inv}\subset B_{rt}.
\end{gather}
For \eqref{multiplierproduct11} choose nets $\{x_i\}$ in $B_r$ and
$\{y_j\}$ in $B_t$ converging to $m_r$ and $m_t$ strictly in $M(B_r)$
and $M(B_t)$, respectively, hence weak* in $B^{**}$. Then
\[ m_{rt}=m_rm_t=\text{weak*}\lim_i\lim_j x_iy_j \]
and
$x_iy_j\in B_rB_t\subset B_{rt}\subset p_{rt}B^{**}p_{t\inv r\inv}$,
so \eqref{multiplierproduct11} holds.
For \eqref{multiplierproduct12}, recall from Lemma \ref{module iso} that
(\ref{multiplier condition for s}) implies $B_r^*m_r=D_{r^{-1}}$ and
$D_r=B_rm_r$. Since we have already verified that (\ref{reducedword})
applies, we deduce that
\begin{equation*}
\begin{split}
D_{rt}m_{rt}
& =B_{rt}B^*_{rt}m_rm_t
=B_rB_tB^*_tB^*_rm_rm_t
=B_rD_tD_{r\inv}m_t\\
& =B_rD_{r\inv}D_tm_t
=B_rB_t
=B_{rt},
\end{split}
\end{equation*}
and similarly for $m_{rt}D_{t\inv r\inv}$.

For \eqref{multiplierproduct2}, since $m_r$ and $m_t$ are partial
isometries with commuting domain and range projections, $m_{rt}=m_rm_t$
is a partial isometry.  Since $m_{rt}\in M(B_{rt})$, we have
$m_{rt}m^*_{rt}\leq p_{rt}$.
For the opposite inequality, it suffices to show
$D_{rt}m_rm_tm^*_tm^*_r=D_{rt}$, and again we use (\ref{reducedword}):
\begin{equation*}
\begin{split}
D_{rt}m_rm_tm^*_tm^*_r
& = B_{rt}B^*_{rt}m_rp_tm^*_r
= B_rB_tB^*_tB^*_rm_rp_tm^*_r\\
& = B_rD_tD_{r\inv}p_tm^*_r
= B_rD_{r\inv}D_tp_tm^*_r \\
& = B_rD_{r\inv}D_tm^*_r
= B_rD_tD_{r\inv}m^*_r \\
& = B_rB_tB^*_tB^*_r
= D_{rt}.
\end{split}
\end{equation*}
Thus \eqref{multiplierproduct2} holds, and we have proved
\eqref{multiplier condition for s} for all $s\in\F_n$.

We still need to verify that $m$ is a partial representation.
\eqref{identity} and \eqref{inverse} are obvious, so it remains to show
\eqref{hom} for all $s,t\in\F_n$. Let $s=s_1\cdots s_k$
and $t=t_1\cdots t_j$ be the reduced spellings. Then the reduced
spelling of $st$ is of the form
\[ st=s_1\cdots s_it_{k-i+1}\cdots t_j
\spacetext{for some}i\le k. \]
Let
\[ u=s_1\cdots s_i,\quad v=s_{i+1}\cdots s_k,\spacetext{and}
w=t_{k-i+1}\cdots t_j. \]
Then $s=uv$, $t=v\inv w$, and $st=uw$, with no cancellations, so
\begin{equation*}
\begin{split}
m_sm_tm^*_tm^*_s
& =m_{uv}m_{v\inv w}m^*_{v\inv w}m^*_{uv}
=m_um_vm^*_vm_wm^*_wm_vm^*_vm^*_u\\
& =m_um_vm^*_vm_wm^*_wm^*_u
=m_{uv}m_{v\inv w}m^*_{uw}
=m_sm_tm^*_{st}.\qed
\end{split}
\end{equation*}
\noqed
\end{proof}

\begin{rem}
When $n=1$, the above theorem includes
\cite[Theorem 4.21]{exe:partial}, and our proof in this case is simpler than
his since we
can use the multiplier bimodule to go straight to the partial
isometries $m_s$.
\end{rem}

\section{Applications and examples}

\subsection*{(a) Cuntz algebras}

The \emph{Toeplitz-Cuntz algebra} $\TO_n$ is the universal \cstar
algebra generated by $n$ isometries $s_i$ such that $\sum_i s_is_i^*$
is a proper projection; Cuntz showed that any $n$ isometries $S_i$ on
Hilbert space generate a faithful representation of $\TO_n$ if $\sum_i
S_iS_i^*<1$ \cite{cun}. The \emph{Cuntz algebra} $\O_n$ is similarly
generated by any family $\{S_i\}$ of isometries satisfying $\sum_i
S_iS_i^*=1$ \cite{cun}.  If $g_i$ are generators of $\F_n$, then
$s_i\otimes g_i\in \TO_n\otimes C^*(\F_n)$ is also a Toeplitz-Cuntz
family of isometries, and hence there is a faithful, unital homomorphism
$\delta\colon \TO_n\to\TO_n\otimes C^*(\F_n)$ such that
$\delta(s_i)=s_i\otimes g_i$. Since
\[ (i\otimes\delta_{\F_n})\circ\delta(s_i)=
i\otimes\delta_{\F_n}(s_i\otimes g_i)= 
s_i\otimes g_i\otimes g_i=(\delta\otimes i)\circ\delta(s_i), \]
$\delta$ is a coaction of $\F_n$ on $\TO_n$. Since
$\{s_i\otimes\lambda_{g_i}\}$ is a Toeplitz-Cuntz family,
$(i\otimes\lambda)\circ\delta$ is faithful, and
$\delta$ is normal. There is a similar coaction on $\O_n$. We intend to
apply
\thmref{elementaryduality} to these coactions.

We first recall some standard notation and facts about the
Toeplitz-Cuntz family $\{s_i\}$. If
$\mu=(\mu_1,\mu_2,\cdots,\mu_{|\mu|})$ is a multi-index, we write
$s_\mu$ for the isometry $s_\mu=s_{\mu_1}s_{\mu_2}\cdots
s_{\mu_{|\mu|}}$ in $\TO_n$, and $g_\mu$ for the word
$g_{\mu_1}g_{\mu_2}\cdots g_{\mu_{|\mu|}}$ in $\F_n$, so that
$\delta(s_\mu)=s_\mu\otimes g_\mu$. If we realise $\TO_n$ on Hilbert
space, the isometries $s_i$ have orthogonal ranges, and hence satisfy
$s_i^*s_j=0$ for $i\neq j$; it follows that every non-zero word in the
$s_i$ and $s_j^*$ collapses to one of the form $s_\mu s_\nu^*$, for
which we have $\delta(s_\mu s_\nu^*)=s_\mu s_\nu^*\otimes g_\mu
g_\nu^{-1}$. A product $(s_\mu s_\nu^*)(s_\alpha s_\beta^*)$ is
non-zero if and only if the multi-indices $\nu$ and $\alpha$ agree as
far as possible; in particular, we have
\[ (s_\mu s_\mu^*)(s_\nu s_\nu^*)
=\begin{cases}
s_\nu s_\nu^*
&\text{if $\nu=(\mu_1,\cdots,\mu_{|\mu|},
\nu_{|\mu|+1},\cdots,\nu_{|\nu|})$};\\
0&\text{if $\nu_i\neq\mu_i$ for some $i\leq\min(|\mu|,|\nu|)$};\\
s_\mu s_\mu^*&
\text{if $\mu=(\nu_1,\cdots,\nu_{|\nu|},
\mu_{|\nu|+1},\cdots,\mu_{|\mu|})$}.
\end{cases} \]
It follows that $D=\clsp\{s_\mu s_\mu^*\}$ is a commutative
\cstar subalgebra of $\TO_n$, which is called the \emph{diagonal subalgebra}.
By convention, we write $s_\emptyset=1$, so that $D$ contains the identity
of $\TO_n$.

\begin{cor}
There is a multiplicative partial action $\alpha$ of $\F_n$ on the
diagonal subalgebra $D$ such that $(\TO_n,\delta)$ is isomorphic to
$(D\times_\alpha\F_n,\hat\alpha)$. Similarly, $\O_n$ is isomorphic
to the partial crossed product of its diagonal subalgebra by a
multiplicative partial action of $\F_n$.
\end{cor}

\begin{proof}
Let $B=\TO_n$.
Since the words $s_\mu s_\nu^*$ span a dense subspace
of $B$, and the projections onto the spectral subspaces $B_s$ are
continuous, the equation $\delta(s_\mu s_\nu^*)=s_\mu s_\nu^*\otimes
g_\mu g_\nu^{-1}$ implies that for each $s\in \F_n$,
\[ B_s=\clsp\{s_\mu s_\nu^*:g_\mu g_\nu^{-1}=s\}. \]
In particular, $B^\delta=B_e=D$, and 
\[ B_{g_i}=\clsp\{s_i s_\mu s_\mu^*\}=s_iD\cong
D=B_{g_i}^*B_{g_i}=D_{{g_i^{-1}}}, \]
as Hilbert $D$-modules. Since the isometries $s_i$ themselves generate
$B$, it follows from \thmref{elementaryduality} that there is a
multiplicative partial action $\alpha$ on $D$ such that
$(B,\delta)\cong(D\times_{\alpha,r}\F_n,\hat\alpha^n)$.

To see that $D\times_\alpha\F_n=D\times_{\alpha,r}\F_n$ in this case,
note that since $D$ is generated by $\{p_s\}_{s\in\F_n}$,
$D\times_\alpha\F_n$ is generated by $\{m_s\}_{s\in\F_n}$. Since the
partial action $\alpha$ is multiplicative, it follows from
\lemref{elementaryprop} that $D\times_\alpha\F_n$ is actually generated
by the Toeplitz-Cuntz family $\{m_{g_1},\ldots,m_{g_n}\}$, and by
Cuntz's Theorem is therefore isomorphic to $\TO_n$. Thus any
representation $\pi\times v$ of $D\times_\alpha\F_n$ in which $\sum
v_iv_i^*\neq 1$ is faithful, including the regular representation whose
image is $D\times_{\alpha,r}\F_n$. The proof for $\O_n$ is similar.
\end{proof}

Much the same arguments show that the Cuntz-Krieger algebras are partial
crossed products; to avoid repetition, we shall merely realize them as
reduced crossed products.

\begin{cor}
If $A$ is a $\{0,1\}$-matrix satisfying condition \textup(I\textup)
of \cite{cun-kri}, then the Cuntz-Krieger algebra $\O_A$ is isomorphic to a
partial crossed product $D\times_{\alpha,r}\F_n$.
\end{cor}

\begin{proof}
Let $n=|A|$, and let $\{s_i\}$ be a family of partial
isometries generating $\O_A$ and satisfying the Cuntz-Krieger relations
\[ s_i^*s_i=\sum_{j=1}^n A(i,j)s_js_j^*. \]
The universal property of $\O_A$ implies that the map $s_i\mapsto s_i\otimes
g_i$ extends to a coaction $\delta\colon\O_A\to\O_A\otimes\F_n$ with
spectral subspaces
\[ B_{g_i}=\overline{{\rm sp}}\{s_i
s_\mu s_\mu^*:A(i,\mu_1)=A(\mu_k,\mu_{k+1})\}, \]
and $D:=B_e=\clsp\{s_\mu s_\mu^*\}$. Since $s_i=\sum_j
s_is_js_j^*$ is in $B_{g_i}$, the spectral subspaces generate $\O_A$. The
ideals $D_{g_i^{-1}}:=B_{g_i}^*B_{g_i}$ are given by
\[ D_{g_i^{-1}}=\clsp\{s_\mu
s_\mu^*:A(i,\mu_1)=A(\mu_k,\mu_{k+1})\}, \]
and $s_i s_\mu s_\mu^*\mapsto s_i^*s_i s_\mu s_\mu^*=s_\mu s_\mu^*$ is a
right Hilbert $D_{g_i^{-1}}$-module isomorphism of $B_{g_i}$ onto
$D_{g_i^{-1}}$. Thus the Corollary follows from
\thmref{elementaryduality}.
\end{proof}

\begin{ex}
We now give some related examples of systems which are not dual to
partial crossed products. First of all, we claim that $\O_n$ is not a
crossed product by a partial action of $\Z$ in such a way that the
gauge action $\alpha$ of $\T$ agrees with the dual action. The spectral
subspaces $B_n$ for the gauge action are given by 
\[ B_n=\clsp\{s_\mu s_\nu^*:|\mu|-|\nu|=n\}, \]
and all the ideals $D_n:=B_nB_n^*$ are equal to the
$AF$-core $B^\alpha$ of $\O_n$. If $(\O_n,\alpha)$ were the dual system of a
partial crossed product $D\times_{\alpha,r}\Z$, then the spectral subspaces
would be isomorphic to $B^\alpha$ as Hilbert $B^\alpha$-modules; but
\[ B_1=\clsp\{s_\mu s_\nu^*:|\mu|-|\nu|=1\}=
\clsp\{s_it:t\in B^\alpha\} \]
is mapped isomorphically to the Hil\-bert $B^\alpha$-mod\-ule $(B^\alpha)^n$ 
via the
map $r\mapsto (s_1^*r,\dots,s_n^*r)$, which has inverse
$(t_1,\dots,t_n)\mapsto\sum s_i t_i$. That $\O_n$ is not a partial
crossed product by $\Z$ underlines that stabilisation is an essential
ingredient in  the comment at the top of
\cite[page 4]{exe:partial} to the effect that the crossed products by
endomorphisms studied by Paschke \cite{pas} fit the mould of
\cite{exe:partial}. (Cuntz's description of $\O_n$ as a crossed product of the
UHF-core $A$ by an endomorphism does not obviously fit the pattern
because the range of the endomorphism is not an ideal in the simple
algebra $A$ (see \cite[Example 3.1]{boy-kes-rae}).

More generally, the coaction of $\F_n$ on $\O_n$ induces a coaction of
any quotient $G$ of $\F_n$, but arguments like those in the previous
paragraph show that these are not typically the dual coaction on some
decomposition $\O_n\cong A\times_{\alpha,r}G$ as a partial crossed
product. For example, if $q\colon\F_3\to \F_2=\langle b_1,b_2\rangle$ is
the quotient map which identifies the first two generators (say
$q(g_1)=q(g_2)=b_1$,
$q(g_3)=b_2$), then the composition $(i\otimes q)\circ\delta:\O_3 \to
\O_3\otimes C^*(\F_2)$ has the Hilbert $B_e$-module
$B_{b_1}=\clsp\{s_it:t\in B_e,i=1,2\}$ isomorphic to
$B_e^2$ rather than $B_e$.
\end{ex}

While it is not an immediate application of our earlier results, it is
interesting to note that similar ideas give a complete characterisation of
the Cuntz algebras in terms of the canonical coaction:

\begin{prop}
Suppose that $B$ is a \cstar algebra with identity,
carrying a coaction $\delta\colon B\to B\otimes C^*(\F_n)$ of $\F_n=\langle
g_1,\cdots,g_n\rangle$. Assume that\textup:
\begin{enumerate}
\item the spectral subspaces $B_{g_i}$ are isomorphic to $B^\delta$ as
right Hilbert $B^\delta$-modules\textup;
\item $B_{g_i}^*B_{g_j}=0$ for $i\neq j$\textup;
\item the spaces $B_{g_i}$ generate $B$ as a \cstar algebra.
\end{enumerate}
Then for every word $\mu$ in the semigroup generated by $\{g_i\}$, $B_\mu
B_\mu^*$ is an ideal in $B^\delta$; as a \cstar algebra, each $B_\mu
B_\mu^*$ has an identity which is a projection $p_\mu$ in $B^\delta$.
\textup(We understand $p_\emptyset$ to be the identity of $B^\delta$.\textup)
Assume further that\textup:
\begin{enumerate}
\setcounter{enumi}{3}
\item $B^\delta=\clsp\{p_\mu:
\mu\text{ is a word in the }g_i\}$.
\end{enumerate}
Then $B$ is isomorphic to $\TO_n$ or $\O_n$. 
\end{prop}

\begin{proof}We begin by observing that the first two assumptions give a
copy of $\O_n$ or $\TO_n$ inside $B$. If $\psi_i\colon B^\delta\to
B_{g_i}$ is the isomorphism guaranteed by (i), then $s_i=\psi_i(1)$
satisfies
\[ s_i^*s_i=\langle s_i,s_i\rangle_{B^\delta}=
\langle \psi_i(1),\psi_i(1)\rangle_{B^\delta}=\langle1,1\rangle_{B^\delta}=1, \]
and hence is an isometry. Since $\psi_i(b)=\psi_i(1b)=s_ib$, we have
$B_{g_i}=s_iB^\delta$. Thus (ii) forces $s_i^*s_j=0$, and $\{s_i\}$ is a
Toeplitz-Cuntz family.

We next claim that $B_\mu=B_{\mu_1}\cdots B_{\mu_{|\mu|}}=s_\mu B^\delta$ for
all words
$\mu$ in $\{g_i\}$. For any $s,t\in \F_n$, we have $B_sB_t\subset B_{st}$,
so the problem is to show that $B_{\mu\nu}\subset B_\mu B_\nu$ for all
words $\mu,\nu$ in $\{g_i\}$. Note that $B_\nu^*B_\nu=B^\delta$, because
$s_\nu^*s_\nu\in B_\nu^*B_\nu$ and $B_\nu^*B_\nu$ is an ideal. Thus 
\[ B_{\mu\nu}=B_{\mu\nu}B_\nu^*B_\nu=B_{\mu\nu}B_{\nu^{-1}}B_\nu\subset
B_{\mu\nu{\nu^{-1}}}B_\nu=B_\mu B_\nu, \]
establishing the first equality. For the second, note that we certainly
have $s_\mu B^\delta\subset B_\mu$. To see that $s_\mu B^\delta$ is all of
$B^\delta$, note that for any $i$, $B^\delta s_i$ is contained in the
spectral subspace $B_{g_i}$, which we know is $s_i B^\delta$. Thus, from
the first equality, we have
\begin{align*}
B_\mu
& =B_{\mu_1}\cdots B_{\mu_{|\mu|}}
=s_{\mu_1}B^\delta s_{\mu_2}B^\delta\cdots s_{\mu_{|\mu|}}B^\delta\\
&\subset s_{\mu_1}s_{\mu_2}\cdots s_{\mu_{|\mu|}}B^\delta
=s_\mu B^\delta,
\end{align*}
justifying the claim.

It follows from the claim that $B_\mu B_\mu^*=s_\mu B^\delta s_\mu^*$,
which has identity $s_\mu s_\mu^*$, and hence (iv) says precisely that
$B^\delta=\overline{{\rm sp}}\{s_\mu s_\mu^*\}$. Thus (iii) implies
that the isometries $s_i$ generate $B$, and $B$ is either $\TO_n$ or
$\O_n$ depending on whether $\sum s_i s_i^*<1$ or $\sum s_is_i^*=1$.
\end{proof}

\subsection*{(b) Wiener-Hopf \cstar algebras}

We now consider the quasi-lattice ordered groups $(G,P)$ of Nica
\cite{nic}. Thus $P$ is a subsemigroup of a discrete group $G$ such
that $P\cap P^{-1}=\{e\}$, and the (right) order on on $G$ defined by
$s\leq t\iff s^{-1}t\in P$ has the following property: if
$s_1,s_2,\ldots,s_n$ have a common upper bound in $P$, they also have a
least upper bound $s_1\vee s_2\vee\cdots\vee s_n$ in $P$. The
individual elements of $G$ which have upper bounds in $P$ are precisely
those in $PP^{-1}=\{pq^{-1}:p,q\in P\}$, and we follow Nica in writing
$\sigma(s)$ for the least upper bound in $P$ of $s\in PP^{-1}$, and $\tau(s)$
for $s^{-1}\sigma(s)$. In general, there are many ways of writing a
given element of $PP^{-1}$ ($pq^{-1}=pr(qr)^{-1}$ for any $r$), and one
should think of $s=\sigma(s)\tau(s)^{-1}$ as the most efficient. We
refer to \cite[Section 2]{nic} for the basic properties and examples.

The \emph{Wiener-Hopf $C^*$-algebra} of a quasi-lat\-tice
ordered group is the $C^*$-algebra $\W(G,P)$ of operators on $l^2(P)$ 
generated by the isometries $\{W_p:p\in P\}$, where
\[ (W_p\xi)(q)=\begin{cases}
\xi(p^{-1}q)&\text{if $p^{-1}q\in P$,}\\
0&\text{otherwise}.
\end{cases} \]
It turns out that the family $\{W_pW_q^*:p,q\in P\}$ spans a dense
subspace of $\W(G,P)$ \cite[Proposition 3.2]{nic}. The \emph{diagonal
subalgebra} is $\D=\clsp\{W_pW_p^*\}$ (see \cite[Section 3]{nic}).

\begin{prop}
\label{reducedcoaction}
There is a normal coaction $\delta$ of $G$ on $\W(G,P)$ such that
\begin{equation}
\label{coactioneq}
\delta(W_pW_q^*)=W_pW_q^*\otimes pq^{-1}\spacetext{for}p,q\in P.
\end{equation}
\end{prop}

\begin{proof}By \cite[Theorem 4.7]{qui:fullred}, it suffices to show there
is a reduced coaction $\delta^r$ of $G$ on $\W(G,P)$ such that
\begin{equation}
\label{reducedcoactioneq}
\delta^r(W_pW_q^*)=W_pW_q^*\otimes\lambda_{pq^{-1}}\spacetext{for}p,q\in P.
\end{equation}
Since $\W(G,P)$ is by definition a subalgebra of $B(l^2(P))$,
the minimal tensor product $\W(G,P)\otimes C^*_r(G)$ by definition acts
on $l^2(P)\otimes l^2(G)=l^2(P\times G)$. We define an operator $W_P$
on $l^2(P\times G)$ by $(W_P\xi)(p,s)=\xi(p,p^{-1}s)$; $W_P$ is unitary
with $W_P^*=W_P^{-1}$ given by $(W_P^*\xi)(p)=\xi(p,ps)$. An easy
calculation shows that
\begin{equation}
\label{conjugation}
\begin{split}
\big(W_P(W_p\otimes1)W_P^*\xi\big)(q,s)& =
\begin{cases}
\xi(p^{-1}q,p^{-1}s)&\text{if $p^{-1}q\in P$,}\\
0&\text{otherwise}\\
\end{cases}\\
& =(W_p\otimes\lambda_p)(\xi)(q,s).
\end{split}
\end{equation}
Since the elements $W_pW_q^*$ span a dense subspace of $\W(G,P)$, the
isometric map $T\mapsto W_P(T\otimes 1)W_P^*$ extends to a unital
homomorphism $\delta^r\colon\W(G,P)\to B(l^2(P\times G))$ with range in
$\W(G,P)\otimes C^*_r(G)$, and \eqref{conjugation} implies
\eqref{reducedcoactioneq}. The coaction identity $(\delta^r\otimes
i)\circ\delta^r=(i\otimes\delta_G^r)\circ\delta^r$ follows easily from
\eqref{reducedcoactioneq}.
\end{proof}

\begin{thm}
\label{weinerhopfduality}
Let $\W(G,P)=C^*(W_p:p\in P)$ be the Wiener-Hopf $C^*$-algebra of a
quasi-lattice ordered group.  Then there is a partial
action $\alpha$ of $G$ on the diagonal subalgebra $\D$ such that the
cosystem $(\W(G,P),G,\delta)$ of Proposition
\textup{\ref{reducedcoaction}} is isomorphic to
$(\D\times_{\alpha,r}G,G,\hat\alpha^n)$.
\end{thm}

\begin{proof}
Let $B$ denote $\W(G,P)$. We aim to apply \thmref{duality}, so we need
a partial representation $m$ of $G$ in $B^{**}$ satisfying
\eqref{multiplier condition}.  For this, we need to identify the
spectral subspaces $B_s$.

\begin{lem}
The spectral subspaces of $\delta$ are given by
\[ B_s=\begin{cases}
W_{\sigma(s)}\D W^*_{\tau(s)}&\text{if }s\in PP\inv,\\
0&\text{otherwise.}
\end{cases} \]
In particular, the fixed-point algebra $B^\delta=B_e$ is the diagonal
subalgebra $\D$.
\end{lem}

\begin{proof}Because $\W(G,P)=\clsp\{W_pW_q^*\}$, \eqref{coactioneq} and
the continuity of the projection $\delta_s=(i\otimes\chi_s)\circ
\delta$ onto $B_s$ imply
\[ B_s=\begin{cases}
\clsp\{W_pW_q^*:pq^{-1}=s\}&\text{if $s\in PP^{-1}$,}\\
0&\text{otherwise.}
\end{cases} \]
Hence, $B_s$ is precisely the subspace $\D_s$ described in
\cite[Section 3.4]{nic}, so the lemma follows from \cite[Section
3.5]{nic}.
\end{proof}

For $s\in G$ define
\[ m_s=\begin{cases}
W_{\sigma(s)}W^*_{\tau(s)}&\text{if }s\in PP\inv,\\
0&\text{otherwise.}
\end{cases} \]
The above lemma tells us that the ideals $D_s=B_sB^*_s$ of $\D$ are
given by
\[ D_s=\begin{cases}
W_{\sigma(s)}\D W^*_{\sigma(s)}&\text{if }s\in PP\inv,\\
0&\text{otherwise,}
\end{cases} \]
so
\begin{equation*}
\begin{split}
m_sm^*_s
& =\begin{cases}
W_{\sigma(s)}W^*_{\sigma(s)}&\text{if }s\in PP\inv,\\
0&\text{otherwise}
\end{cases}\\
& =p_s.
\end{split}
\end{equation*}
In particular, \eqref{identity} holds. Since
\begin{equation*}
\begin{split}
m_{s\inv}
& =\begin{cases}
W_{\sigma(s\inv)}W^*_{\tau(s\inv)}&\text{if }s\in PP\inv,\\
0&\text{otherwise}
\end{cases}\\
& =\begin{cases}
W_{\tau(s)}W^*_{\sigma(s)}&\text{if }s\in PP\inv,\\
0&\text{otherwise}
\end{cases}\\
& =m^*_s,
\end{split}
\end{equation*}
\eqref{inverse} holds as well. Furthermore, the above formulas imply
\eqref{multiplier condition}. It remains to verify \eqref{hom}: for
$s,t\in G$ we must show $m_sm_t\preceq m_{st}$. We may assume $s,t,s\inv
t\in PP\inv$, since $m_sm_t=0$ otherwise. Then
\begin{equation*}
\begin{split}
m_sm_t
& =W_{\sigma(s)}W^*_{\tau(s)}W_{\sigma(t)}W^*_{\tau(t)}\\
& =W_{\sigma(s)\tau(s)\inv(\tau(s)\vee\sigma(t))}
W^*_{\tau(t)\sigma(t)\inv(\sigma(t)\vee\tau(s))},
\end{split}
\end{equation*}
by \cite[Equation(5)]{nic}, while
\[ m_{st}=W_{\sigma(st)}W^*_{\tau(st)}. \]
Since $W_pW^*_q\preceq W_uW^*_v$ whenever $p,q,u,v\in P$, $pq\inv=uv\inv$,
and $u\le p$, the inequality follows.

Thus \thmref{weinerhopfduality} follows from \thmref{duality}.
\end{proof}

\begin{rem}Nica also associates to each $(G,P)$ a universal
$C^*$-algebra $C^*(G,P)$ whose representations are given by
representations $V$ of $P$ as isometries satisfying the covariance
condition
\[ V_pV_p^*V_qV_q^*=\begin{cases}
V_{p\vee q}V_{p\vee q}^*&\text{if $p\vee q$ exists,}\\
0&\text{otherwise.}
\end{cases} \]
If $V\colon P\to C^*(G,P)$ is the universal such representation, the
map $s\mapsto V_s\otimes s$ is also covariant, and hence there is a
coaction $\delta\colon C^*(G,P)\to C^*(G,P)\otimes C^*(G)$ such that
$\delta(V_s)=V_s\otimes s$. It is not obvious that this coaction will
be normal, and its normalisation could coact on a proper quotient of
$C^*(G,P)$. However, the theory of \cite{nic} suggests that in many
cases of interest the Wiener-Hopf representation induces an isomorphism
of $C^*(G,P)$ onto $\W(G,P)$; a general theorem along these lines is
proved in \cite{lr}, from which Cuntz's Theorem \cite{cun} and other
related results follow.
\end{rem}


\providecommand{\bysame}{\leavevmode\hbox to3em{\hrulefill}\thinspace}

\end{document}